\documentclass[11pt,a4paper]{article}

\usepackage{mathtools}
\usepackage{authblk} 
\usepackage{algpseudocode,algorithmicx,algorithm}

\usepackage{mathrsfs}
\usepackage{latexsym,bm}
\usepackage{amsmath,amsfonts,amsmath,amssymb,amsthm}
\usepackage{extarrows}

\usepackage{graphicx,subfigure,epstopdf,float}
\usepackage{enumerate,cases,multirow}
\usepackage{makecell}
\usepackage{caption}

\usepackage{longtable,colortbl,arydshln,threeparttable}
\definecolor{mygray}{gray}{.9}

\usepackage{indentfirst}
\usepackage[top=25mm,bottom=20mm,left=25mm,right=20mm]{geometry}
\usepackage{cite}

\usepackage{listings}

\usepackage{makeidx}        
\usepackage{booktabs}
\usepackage[bookmarksnumbered,colorlinks,citecolor=red,linkcolor=red,hyperindex,linktocpage=true]{hyperref}

\newcommand{\bb}{\boldsymbol}

\def \d {\mathrm{d}}

\newcounter{parentalgorithm}

\makeatother

\theoremstyle{remark}
\newtheorem{remark}{\bf Remark}[section]

\numberwithin{equation}{section}

\usepackage{algorithm}
\usepackage{listings,xcolor} 
\usepackage{textcomp}
\usepackage[framed,numbered,autolinebreaks,useliterate]{mcode} 

\graphicspath{{images/}}

\begin{document}

\title{\bf Unified Implementation of Finite Element Methods involving Jumps and Averages in Matlab}
\author{Yue Yu\thanks{terenceyuyue@sjtu.edu.cn}}
\affil{School of Mathematical Sciences, Institute of Natural Sciences, MOE-LSC, Shanghai Jiao Tong University, Shanghai, 200240, P. R. China.}
\date{}
\maketitle

\begin{abstract}
  We provide unified implementations of the finite element methods involving jumps and averages in Maltab by combing the use of the software package varFEM, including the adaptive finite element methods for the Poisson equation and the $C^0$ interior penalty methods for the biharmonic equation. The design ideas can be extended to other Galerkin-based methods, for example, the discontinuous Galerkin methods and the virtual element methods.
\end{abstract}

\section{Implementation of finite element methods in varFEM package} \label{sect:varFEM}

In this section we briefly review the finite element method and introduce the numerical implementation in varFEM, a Matlab software package for the finite element method.
For simplicity, we consider the Poisson equation with the homogeneous Dirichlet boundary conditions.

\subsection{The finite element method for the Poisson equation}

Let $\Omega$ be a bounded Lipschitz domain in $\mathbb{R}^2$ with polygonal boundary $\partial \Omega$. Consider the Poisson equation
\begin{equation}\label{modelP}
\begin{cases}
   - \Delta u = f \quad & \text{in}~~\Omega ,  \\
  u = 0\quad & \text{on}~~\partial \Omega,
\end{cases}
\end{equation}
where $f\in L^2(\Omega)$ is a given function. The continuous variational problem is to find $u \in V: = H_0^1(\Omega )$ such that
\[a(u,v) = \ell(v),\quad v \in V,\]
where
\[a(u,v) = \int_\Omega  \nabla u \cdot \nabla v {\rm d}\sigma, \qquad \ell (v) = \int_\Omega  fv {\rm d}\sigma.\]

For the finite element discretization, we discuss the conforming Lagrange elements. Let $\mathcal{T}_h$ be a shape regular triangulation. The generic element will be denoted as $K$ in the sequel. Define
\[V_h = \{ v \in V: v |_K \in \mathbb{P}_k(K), \quad K \in \mathcal{T}_h \},\]
where $k \le 3$. The discrete problem is to seek $u_h \in V_h$ satisfying
\begin{equation}\label{FEM}
a(u_h,v) = \ell (v), \qquad v \in V_h.
\end{equation}

\subsection{Data structures for triangular meshes} \label{subsec:datastructure}

\begin{figure}[H]
  \centering
  \includegraphics[scale = 0.45]{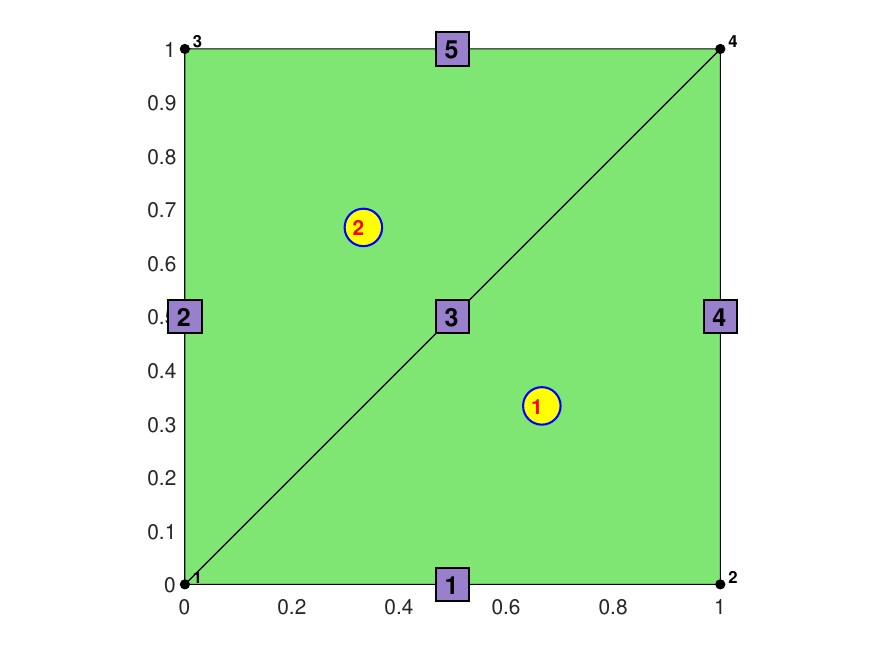}\\
  \caption{Illustration of the data structures}\label{fig:datastructure}
\end{figure}

We adopt the data structures given in $i$FEM \cite{ChenL-iFEM-2009}. All related data are stored in the Matlab structure \mcode{Th}, which is computed by using the subroutine \mcode{FeMesh2d.m} as
\vspace{-0.8cm}
\begin{lstlisting}
Th = FeMesh2d(node,elem,bdStr);
\end{lstlisting}
where the basic data structures \mcode{node} and \mcode{elem} are generated by
\vspace{-0.8cm}
\begin{lstlisting}
[node,elem] = squaremesh([x1 x2 y1 y2], h1, h2);
\end{lstlisting}
For clarity, we take a simple mesh shown in Fig.~\ref{fig:datastructure} as an example.

The triangular meshes are represented by two basic data structures \mcode{node} and \mcode{elem}, where \mcode{node} is an $\mcode{N} \times 2$ matrix with the first and second columns contain $x$- and $y$-coordinates of the nodes in the mesh, and \mcode{elem} is an $\mcode{NT} \times 3$ matrix recording the vertex indices of each element in a counterclockwise order, where \mcode{N} and \mcode{NT} are the numbers of the vertices and triangular elements. For the mesh given in Fig.~\ref{fig:datastructure},
\[\mcode{elem} = \begin{bmatrix} 2 & 4 & 1\\ 3 & 1 & 4 \end{bmatrix}\]

In the current version, we only consider the $\mathbb{P}_k$-Lagrange finite element spaces with $k$ up to 3. In this case, there are two important data structures \mcode{edge} and \mcode{elem2edge}. In the matrix \mcode{edge(1:NE,1:2)}, the first and second rows contain indices of the starting and ending points. The column is sorted in the way that for the $k$-th edge, \mcode{edge(k,1)<edge(k,2)} for $k=1,2,\cdots,\mcode{NE}$. For the given triangulation,
\[\mcode{edge} = \begin{bmatrix}
1	&   2 \\
1   &	3 \\
1   &	4 \\
2   &	4 \\
3   &	4
\end{bmatrix}.\]
The matrix \mcode{elem2edge} establishes the map of local index of edges in each triangle to its global index in matrix \mcode{edge}.
By convention, we label three edges of a triangle such that the $i$-th edge is opposite to the $i$-th vertex. For the given mesh,
\[\mcode{elem2edge} = \begin{bmatrix}
3   &	1 &	4 \\
3   &	5 &	2 \\
\end{bmatrix}.\]
In some cases, we may need to specify the left and right triangles sharing the same edge. For this purpose, we introduce another data structure \mcode{edge2elem}, which is an $\mcode{NE} \times  4$ matrix such that \mcode{edge2elem(k,1)} and \mcode{edge2elem(k,2)} are two triangles sharing the $k$-th edge for an interior edge. If the $k$-th edge is on the boundary, then we set \mcode{edge2elem(k,1) = edge2elem(k,2)}. For convenience, we also record the local index of the edge on the left and right triangles in \mcode{edge2elem(k,3)} and \mcode{edge2elem(k,4)}, respectively. For the given mesh,
\[\mcode{edge2elem} = \begin{bmatrix}
1 &	1 & 2 & 2 \\
2 & 2 & 3 & 3 \\
1 & 2 & 1 & 1 \\
1 & 1 & 3 & 3 \\
2 & 2 & 2 & 2 \\
\end{bmatrix}.\]
We refer the reader to \url{https://www.math.uci.edu/~chenlong/ifemdoc/mesh/auxstructuredoc.html} for some detailed information.

To deal with boundary integrals, we first extract the boundary edges from \mcode{edge} and store them in matrix \mcode{bdEdge}.
In the input of \mcode{FeMesh2d.m}, the string \mcode{bdStr} is used to indicate the interested boundary part in \mcode{bdEdge}. For example, for the unit square $\Omega = (0,1)^2$,
\begin{itemize}
  \item \mcode{bdStr = 'x==1'} divides \mcode{bdEdge} into two parts: \mcode{bdEdgeType\{1\}}
gives the boundary edges on $x=1$, and \mcode{bdEdgeType\{2\}} stores the remaining part.
  \item \mcode{bdStr = \{'x==1','y==0'\} } separates the boundary data structure \mcode{bdEdge} into three parts: \mcode{bdEdgeType\{1\}} and \mcode{bdEdgeType\{2\}} give the boundary edges on $x=1$ and $y=0$, respectively, and \mcode{bdEdgeType\{3\}} stores the remaining part.
  \item \mcode{bdStr = []} implies that \mcode{bdEdgeType\{1\} = bdEdge}.
\end{itemize}

We also use \mcode{bdEdgeIdxType} to record the indices in matrix \mcode{edge}, and \mcode{bdNodeIdxType} to store the node indices for respective boundary parts. Note that we determine the boundary of interest by the coordinates of the midpoint of the edge, so \mcode{'x==1'} can also be replaced by a statement like \mcode{'x>0.99'}.

\subsection{Implementation of the FEM in varFEM} \label{subsect:implementationFEM}

FreeFEM is a popular 2D and 3D partial differential equations (PDE) solver based on finite element methods \cite{FreeFEM}, which has been used by thousands of researchers across the world. The highlight is that the programming language is consistent with the variational formulation of the underlying PDEs, referred to as the variational formulation based programming in \cite{varFEM}, where we have developed an FEM package in a similar way of FreeFEM using the language of Matlab, named varFEM. The similarity here only refers to the programming style of the main or test script, not to the internal architecture of the software.

Consider the Poisson equation with the Dirichlet boundary condition on the unit square. The exact solution is given by
\[u(x,y) = xy(1 - x)(1 - y){\text{exp}}\left(  - 1000((x - 0.5)^2 + (y - 0.117)^2) \right).\]
The PDE data is generated by \mcode{pde = Poissondata\_afem()}.  The function file is simply given as follows.
\vspace{-0.8cm}
\begin{lstlisting}
function uh = varPoisson(Th,pde,Vh,quadOrder)

%% Assemble stiffness matrix
Coef  = 1;
Test  = 'v.grad';
Trial = 'u.grad';
kk = assem2d(Th,Coef,Test,Trial,Vh,quadOrder);

%% Assemble right hand side
Coef = pde.f;  Test = 'v.val';
ff = assem2d(Th,Coef,Test,[],Vh,quadOrder);

%% Apply Dirichlet boundary conditions
g_D = pde.g_D;
on = 1;
uh = apply2d(on,Th,kk,ff,Vh,g_D);
\end{lstlisting}

In the above code, the structure \mcode{pde} stores the information of the PDE, including the exact solution \mcode{pde.uexact}, the gradient \mcode{pde.Du}, etc.  We set up the triple \mcode{(Coef,Test,Trial)} for the coefficients, test functions and trial functions in variational form, respectively. It is obvious that \mcode{v.grad} is for $\nabla v$ and \mcode{v.val} is for $v$ itself. The routine \mcode{assem2d.m}  computes the stiffness matrix corresponding to the bilinear form on the two-dimensional region, i.e.
\[A = (a_{ij}),\quad a_{ij} = a(\Phi _j,\Phi _i),\]
where $\Phi _i$ are the global shape functions of the finite element space \mcode{Vh}. The integral of the bilinear form, as
$(\nabla \Phi _i, \nabla \Phi _j)_\Omega$, is approximated by using the Gaussian quadrature formula with \mcode{quadOrder} being the order of accuracy.

We remark that the \mcode{Coef} has three forms:
\begin{enumerate}
  \item A function handle or a constant.
  \item The numerical degrees of freedom of a finite element function.
  \item A coefficient matrix \mcode{CoefMat} resulting from the numerical integration.
\end{enumerate}
In the computation, the first two forms in fact will be transformed to the third one. Given a function $c(p)$, where $p = (x,y)$, the coefficient matrix is in the following form
\begin{equation}\label{coefmat}
\mcode{CoefMat} = \begin{bmatrix}
c(p_1^1) & c(p_2^1) & \cdots  & c(p_{n_g}^1) \\
c(p_1^2) & c(p_2^2) & \cdots  & c(p_{n_g}^2) \\
\vdots   & \vdots   & \vdots  & \vdots \\
c(p_1^{\text{NT}}) & c(p_2^{\text{NT}}) & \cdots  & c(p_{n_g}^{\text{NT}}) \\
\end{bmatrix}.
\end{equation}
Here, $p_1^i, p_2^i, \cdots, p_{n_g}^i$ are the quadrature points on the element $K_i$.

\begin{figure}[!htb]
  \centering
  \includegraphics[scale=0.45]{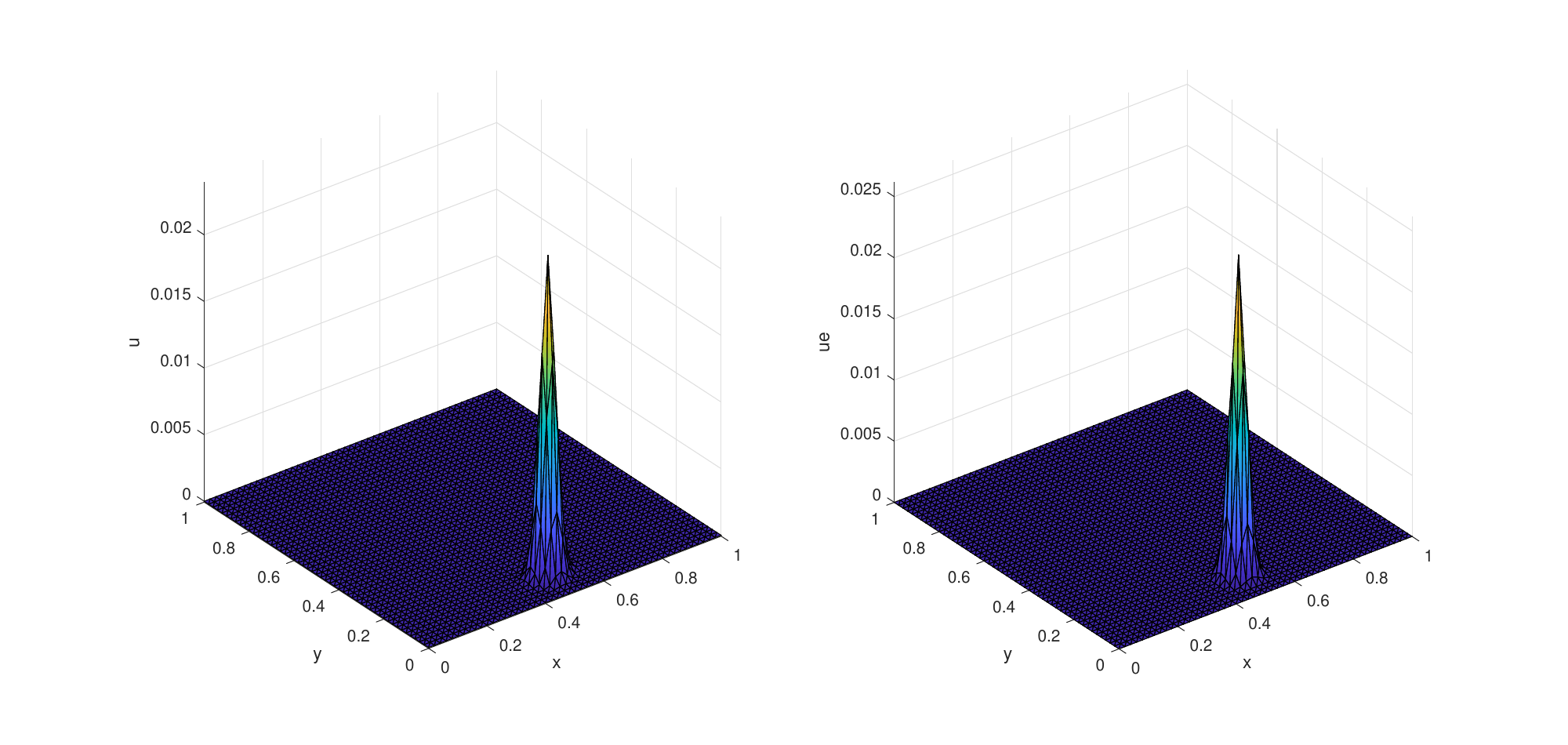}\\
  \caption{Numerical and exact solutions with $h_1 = h_2 = 1/50$}\label{fig:AFEM0}
\end{figure}

We display the numerical result in Fig.~\ref{fig:AFEM0} for the uniform triangular mesh with $h_1 = h_2 = 1/50$ generated by
\vspace{-0.8cm}
\begin{lstlisting}
[node,elem] = squaremesh([0 1 0 1], 1/50, 1/50);
\end{lstlisting}

\subsection{Review on writing the subroutine \mcode{assem2d.m}}

To make the discussions in the subsequent sections more clear, we now introduce the details of writing the subroutine \mcode{assem2d.m} to assemble a two-dimensional scalar bilinear form
\[a(v,u),\quad v = \varphi _i, ~~u = \phi _j,\quad i = 1, \cdots ,m,~~j = 1, \cdots ,n,\]
where the test function $v$ and the trial function $u$ are allowed to match different finite element spaces, which, for example, can be found in mixed finite element methods for Stokes problems. To handle different spaces, we write \mcode{Vh = \{'P1', 'P2' \}} for the input of \mcode{assem2d.m}, where \mcode{Vh\{1\}} is for $v$ and \mcode{Vh\{2\}} is for $u$. For simplicity, it is also allowed to write \mcode{Vh = 'P1'} when $v$ and $u$ are in the same space.

Let us discuss the case where $v$ and $u$ are in the same space. Suppose that the bilinear form contains only first-order derivatives. Then the possible combinations are
\[\int_K {avu} {\rm d}x,\quad \int_K av_xu {\rm d}x,\quad \int_K av_yu {\rm d}\sigma,\]
\[\int_K avu_x {\rm d}\sigma,\quad \int_K av_xu_x {\rm d}\sigma,\quad \int_K av_yu_x {\rm d}\sigma,\]
\[\int_K avu_y {\rm d}\sigma,\quad \int_K av_xu_y {\rm d}\sigma,\quad \int_K av_yu_y {\rm d}\sigma.\]
Of course, we often encounter the gradient form
\[\int_K  a\nabla v \cdot \nabla u {\rm d}\sigma = \int_K  a(v_xu_x + v_yu_y) {\rm d}\sigma.\]

We take the second bilinear form as an example. Let
\[a_K(v,u) = \int_K a v_xu {\rm d}\sigma,\]
and consider the $\mathbb{P}_1$-Lagrange finite element. Denote the local basis functions to be $\phi _1,\phi _2,\phi _3$. Then the local stiffness matrix is
\[A_K = \int_K a  \begin{bmatrix}
  \partial _x\phi _1 \\
  \partial _x\phi _2 \\
  \partial _x\phi _3
\end{bmatrix} \begin{bmatrix}\phi _1 & \phi _2 & \phi _3 \end{bmatrix} {\rm d}\sigma.\]
Let
\[v_1 = \partial _x\phi _1, ~~v_2 = \partial _x\phi _2,~~v_3 = \partial _x\phi _3;\quad
u_1 = \phi _1,~~u_2 = \phi _2,~~u_3 = \phi _3.\]
Then
\[A_K = (k_{ij})_{3 \times 3}, \quad k_{ij} = \int_K av_iu_j {\rm d}\sigma.\]
The integral will be approximated by the Gaussian quadrature rule:
\[
k_{ij} = \int_Kav_iu_j {\rm d}\sigma = | K |\sum\limits_{p = 1}^{n_g} w_pa(x_p,y_p)v_i(x_p,y_p)u_j(x_p,y_p) ,
\]
where $(x_p,y_p)$ is the $p$-th quadrature point. In the implementation, we in advance store the quadrature weights and the values of basis functions or their derivatives in the following form:
\[w_p,\qquad \mcode{vi}(:,p) = \begin{bmatrix}
  v_i|_{(x_p^1,y_p^1)} \\
  v_i|_{(x_p^2,y_p^2)} \\
   \vdots  \\
  v_i|_{(x_p^{\text{NT}},y_p^{\text{NT}})}
\end{bmatrix},\qquad \mcode{uj}(:,p) =  \begin{bmatrix}
  u_j|_{(x_p^1,y_p^1)} \\
  u_j|_{(x_p^2,y_p^2)} \\
   \vdots  \\
  u_j|_{(x_p^{\text{NT}},y_p^{\text{NT}})}
\end{bmatrix},\qquad p = 1, \cdots ,n_g,\]
where $\mcode{vi}$ associated with $v_i$ is of size $\mcode{NT} \times \mcode{ng}$ with the $p$-th column given by $\mcode{vi}(:,p)$.
Let $\mcode{weight} = [w_1, w_2, \cdots, w_{n_g}]$ and \mcode{ww = repmat(weight, NT, 1)}. Then $k_{ij}$ for $a=1$ can be computed as
\vspace{-0.8cm}
\begin{lstlisting}
k11 = sum(ww.*v1.*u1,2);
k12 = sum(ww.*v1.*u2,2);
k13 = sum(ww.*v1.*u3,2);
k21 = sum(ww.*v2.*u1,2);
k22 = sum(ww.*v2.*u2,2);
k23 = sum(ww.*v2.*u3,2);
k31 = sum(ww.*v3.*u1,2);
k32 = sum(ww.*v3.*u2,2);
k33 = sum(ww.*v3.*u3,2);
K = [k11,k12,k13, k21,k22,k23, k31,k32,k33];
\end{lstlisting}
Here we have stored the local stiffness matrix $A_K$ in the form of $[k_{11},k_{12},k_{13},k_{21},k_{22},k_{23},k_{31},k_{32},k_{33}]$, and stacked the results of all cells together. By adding the contribution of the area, one has
\vspace{-0.8cm}
\begin{lstlisting}
Ndof = 3;
K = repmat(area,1,Ndof^2).*K;
\end{lstlisting}

For the variable coefficient case, such as $a(x,y) = x+y$,  one can further introduce the coefficient matrix as
\vspace{-0.8cm}
\begin{lstlisting}
cf = @(pz) pz(:,1) + pz(:,2); % x+y;
cc = zeros(NT,ng);
for p = 1:ng
    pz = lambda(p,1)*z1 + lambda(p,2)*z2 + lambda(p,3)*z3;
    cc(:,p) = cf(pz);
end
\end{lstlisting}
where \mcode{pz} are the quadrature points on all elements. The above procedure can be implemented as follows.
\vspace{-0.8cm}
\begin{lstlisting}
K = zeros(NT,Ndof^2);
s = 1;
v = {v1,v2,v3}; u = {u1,u2,u3};
for i = 1:Ndof
    for j = 1:Ndof
        vi = v{i}; uj = u{j};
        K(:,s) =  area.*sum(ww.*cc.*vi.*uj,2);
        s = s+1;
    end
end
\end{lstlisting}
The bilinear form is assembled by using the built-in function \mcode{sparse.m} as in $i$FEM. In this case, the code is given as
\vspace{-0.8cm}
\begin{lstlisting}
ss = K(:);
kk = sparse(ii,jj,ss,NNdof,NNdof);
\end{lstlisting}
Here, \mcode{(ii, jj)} is called the sparse index.

\begin{remark}\label{rem:Base2d}
In varFEM, we use \mcode{Base2d.m} to load the information of \mcode{vi} and \mcode{uj}, for example, the following code gives the values of $\partial_x \phi$, where $\phi$ is a local basis function.
\vspace{-0.8cm}
\begin{lstlisting}
v = 'v.dx';
vbase = Base2d(v,node,elem,Vh{1},quadOrder);
\end{lstlisting}
\end{remark}

\begin{remark}\label{rem:elem2dof}
The sparse index \mcode{(ii, jj)} can be simply given by the connectivity list \mcode{elem2dof} as
\vspace{-0.8cm}
\begin{lstlisting}
% assembly index
ii = reshape(repmat(elem2dof, Ndof,1), [], 1);
jj = repmat(elem2dof(:), Ndof, 1);
\end{lstlisting}
where \mcode{elem2dof} is an $\mcode{NT} \times \mcode{Ndof}$ matrix, with the $i$-th row representing the index vector of the degrees of freedom. For example, \mcode{elem2dof = elem} for the $\mathbb{P}_1$-Lagrange element.
\end{remark}

\section{Adaptive finite element methods for the Poisson equation} \label{sect:varFEM}

In this section we briefly introduce the ingredients of the adaptive finite element method, and provide the overall structure of the implementation.

\subsection{The ingredients of the adaptive FEM}

For the conforming FEM, one can establish the following residual based a-posteriori error estimate
\[\| u - u_h \|_1 \lesssim \eta (u_h),\]
where $\eta  = \Big( \sum\limits_{K \in \mathcal{T}_h} \eta _K^2  \Big)^{1/2}$ and
\begin{equation}\label{localErrorIndicator}
\eta _K^2 = h_K^2\| f + \Delta u_h\|_{0,K}^2 + \sum\limits_{e \subset \partial K} h_e\| [\partial _{n_e}u_h] \|_{0,e}^2
\end{equation}
are referred to as the global and local error indicators, respectively. Here, $f + \Delta u_h$ is the interior residual, and
\[[\partial _{n_e}u_h] = \partial _{n_e}u_h |_{K^-} - \partial _{n_e}u_h |_{K^+}\]
is the jump term, where $n_e$ is the unit outer normal vector of the target element $K = K^-$, and $K^+$ is the neighbouring triangle sharing $e$ as an edge. Note that for an edge $e$ on the domain boundary $\partial \Omega$, we assume $\partial _{n_e}u_h |_{K_2} = 0$. For some problems, we also need to introduce the average: $\{u\} = \frac{1}{2} (u|_{K^-} + u|_{K^+})$ for an interior edge and $\{u\} = u$ for a boundary edge.

Standard adaptive algorithms based on the local mesh refinement can be written as loops of the form
\[{\bf SOLVE} \to {\bf ESTIMATE} \to {\bf MARK} \to {\bf REFINE}.\]
Given an initial subdivision $\mathcal{T}_0$, to get $\mathcal{T}_{k+1}$ from $\mathcal{T}_k$ we
first solve the FEM problem under consideration to get the numerical solution $u_k$ on $\mathcal{T}_k$. The error is then estimated
by using $u_k$, $\mathcal{T}_k$ and the a posteriori error bound $\eta(u_h)$. The local error bound $\eta_K$ is used to mark a subset $\mathcal{M}$
of elements in $\mathcal{T}_k$ for refinement. The marked triangles and possible more neighboring elements are refined
in such a way that the subdivision meets certain conditions, for example, the resulting triangular mesh is still shape
regular. The above procedures are included in the test script with an overview given as follows
\vspace{-0.8cm}
\begin{lstlisting}
for k = 1:maxIt
    % Step 1: SOLVE
    Th = FeMesh2d(node,elem,bdStr);
    uh = varPoisson(Th,pde,Vh,quadOrder);
    % Step 2: ESTIMATE
    eta = Poisson_indicator(Th,uh,pde,Vh,quadOrder);
    % Step 3: MARK
    elemMarked = mark(elem,eta,theta);
    % Step 4: REFINE
    [node,elem] = bisect(node,elem,elemMarked);
end
\end{lstlisting}

Three important modules are involved: the local error indicator in Step 2, the marking algorithm in Step 3 and the local refinement algorithm in Step 4. We employ the D\"{o}rfler marking strategy to select the subset of elements and then use the newest vertex method to refine the mesh. Note that the subroutines \mcode{mark.m} and \mcode{bisect.m} are extracted from $i$FEM \cite{ChenL-iFEM-2009} with some minor modifications. In this article, we are only concerned with the computation of the error indicator. For the later two steps, we refer the reader to \cite{ChenL-iFEM-2009} for details on the implementation.

\subsection{The unified implementation of the error indicator}

\subsubsection{The computation of the elementwise residuals}

With the help of varFEM, the first term $h_K^2\| f + \Delta u_h\|_{0,K}^2$ in \eqref{localErrorIndicator} can be simply computed as
\vspace{-0.8cm}
\begin{lstlisting}
%% elementwise residuals
fc = interp2dMat(pde.f,'.val',Th,Vh,quadOrder);
uxxc = interp2dMat(uh,'.dxx',Th,Vh,quadOrder);
uyyc = interp2dMat(uh,'.dyy',Th,Vh,quadOrder);
Coef = (fc + uxxc + uyyc).^2;
[~,elemIh] = integral2d(Th,Coef,Vh,quadOrder);
elemRes = diameter.^2.*elemIh;
\end{lstlisting}
In the above code, \mcode{elemIh} stores all the local error indicators $[\eta_{K_1}, \cdots, \eta_{K_{\text{NT}}}]^T$; The function \mcode{interp2dMat} is used to generate the coefficient matrix (see \eqref{coefmat}). It is obvious that the coefficient matrix of $(f + \Delta u_h)^2$ is \mcode{(fc+uxxc+uyyc).^2}, where \mcode{fc}, \mcode{uxxc} and \mcode{uyyc} are the coefficient matrices of $f$, $\partial_{xx}u_h$ and $\partial_{yy}u_h$, respectively.

In what follows, we focus on the unified implementation of the jump term or jump integral $\sum_{e \subset \partial K} h_e\| [\partial _{n_e}u_h] \|_{0,e}^2$. We remark that any other type of jump terms can be easily adapted or designed accordingly as will be seen.

\subsubsection{The elementwise interior and exterior indices of the quadrature points}

The integral over $e$ is calculated by using the one-dimensional Gaussian numerical integration formula
\[\int_e f \d s = |e|( w_1 f(p_1) + w_2 f(p_2) + \cdots + w_{n_g} f(p_{n_g})),\]
where $w_i$ and $p_i$ are the quadrature weights and points on $e$, and $n_g$ is the number of the quadrature points. Note that the endpoints of $e$ are not included.

\begin{figure}[H]
  \centering
  \includegraphics[scale=0.5]{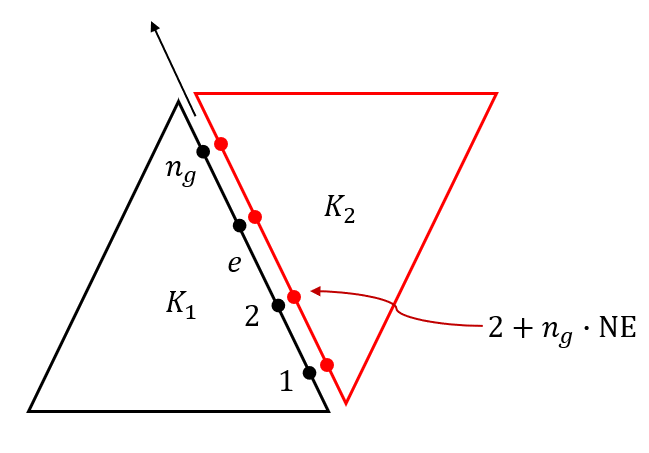}\\
  \caption{Illustration of the quadrature points}\label{fig:quadpoints}
\end{figure}

\begin{enumerate}
  \item Indexing rule for the quadrature points. In Fig.~\ref{fig:quadpoints} the direction of the edge $e$ is specified by the arrow. In the implementation, the direction is determined by the data structure \mcode{edge} which satisfies \mcode{edge(k,1)<edge(k,2)}. For the first edge $e = e_1$ in \mcode{edge}, the quadrature points will be numbered by $1,2,\cdots, n_g$ or $1:n_g$ for short when restricted to the left triangle. Similarly, the quadrature points on the $i$-th edge $e_i$ will be numbered by
  \[p_i^L = (1:n_g) + (i-1)n_g , \qquad i = 1,2,\cdots, \mcode{NE}.\]
  The right-hand restrictions are accordingly numbered as
  \[p_i^R = p_i^L + n_g\cdot \mcode{NE}, \qquad i = 1,2,\cdots, \mcode{NE}.\]
  \item Elementwise sign matrix. To characterize the direction of an edge on an element $K$, one can introduce the elementwise sign matrix
  \vspace{-0.8cm}
  \begin{lstlisting}
% sign of elementwise edges
sgnelem = sign([elem(:,3)-elem(:,2), elem(:,1)-elem(:,3), elem(:,2)-elem(:,1)]);
  \end{lstlisting}
  In some cases, it is better to restore the positive sign for edges on the boundary of the domain, or one can set the signs to be zero for later use:
\vspace{-0.8cm}
\begin{lstlisting}
E = false(NE,1); E(bdEdgeIdx) = 1; sgnbd = E(elem2edge);
sgnelem(sgnbd) = 0;
  \end{lstlisting}
  Here, the data structures \mcode{elem2edge} and \mcode{bdEdgeIdx} have been introduced in Subsect.~\ref{subsec:datastructure}.

  \item Elementwise interior indices. According to the indexing rule, one easily obtains the interior indices of the first sides of all triangles:
    \vspace{-0.8cm}
  \begin{lstlisting}
% first side
e1 = elem2edge(:,1); sgn1 = sgnelem(:,1);
id = repmat(1:ng, NT,1);
id(sgn1<0,:) = repmat((ng:-1:1)+NE*ng, sum(sgn1<0), 1);
elemQuade1 = id + (e1-1)*ng;
  \end{lstlisting}
  Note that for the edges with positive and zero signs the natural indices are $1:n_g$, while for those with negative signs the natural indices are reversed. One can similarly introduce the interior indices of other sides, and hence give the elementwise interior indices
   \vspace{-0.8cm}
  \begin{lstlisting}
elemQuadM = [elemQuade1,elemQuade2,elemQuade3];
  \end{lstlisting}
  where the letter \mcode{M} stands for ``Minus''.

  \item Elementwise exterior indices. To compute the jump, we also need to introduce the elementwise exterior indices \mcode{elemQuadP}, where \mcode{P} is for ``Plus''.  Given an edge, assume that the interior indices of the quadrature points are $i_1,\cdots,i_{n_g}$, and the exterior indices are $i'_1,\cdots,i'_{n_g}$. Since $|i_k-i'_k| = n_g\cdot \mcode{NE}$, one just needs to subtract $n_g\cdot \mcode{NE}$ for those indices in \mcode{elemQuadM} greater than  $n_g\cdot \mcode{NE}$, and add $n_g\cdot \mcode{NE}$ to those indices less than $n_g\cdot \mcode{NE}$.
     \vspace{-0.8cm}
  \begin{lstlisting}
index = ( elemQuadM > ng*NE ) ;
elemQuadP = elemQuadM + (-ng*NE)*index + ng*NE*(~index);
  \end{lstlisting}
  Obviously, for edges on the domain boundary, the exterior indices are greater than $n_g \cdot \mcode{NE}$.
\end{enumerate}

\subsubsection{The elementwise interior and exterior evaluations} \label{subsubsect:MPevaluation}

Given $u_h$, we now determine the interior evaluations \mcode{elemuhM} and the exterior evaluations \mcode{elemuhP} of $u_h$ corresponding to \mcode{elemQuadM} and \mcode{elemQuadP}. To this end, one can compute the evaluations of the basis functions at the quadrature points along the boundary $\partial K$.
Given a triangle $K$, let $\lambda_1$, $\lambda_2$ and $\lambda_3$ be the barycentric coordinate functions. Since $\lambda_i$ are usually used to construct the basis functions for the FEMs, one can first specify the evaluations of $\lambda_i$ and the derivatives $\partial_x \lambda_i$ and $\partial_y \lambda_i$ at the quadrature points.

Given some points $p_1,p_2,\cdots,p_{n_G}$ on the triangle $K$, for example, the quadrature points for the integration of the bilinear forms. In varFEM, the evaluations of $\lambda_i$ are stored in the following form
\[ \begin{bmatrix}
\lambda_1(p_1) & \lambda_2(p_1) & \lambda_3(p_1) \\
\lambda_1(p_2) & \lambda_2(p_2) & \lambda_3(p_2) \\
\vdots         & \vdots         & \vdots         \\
\lambda_1(p_{n_G}) & \lambda_2(p_{n_G}) & \lambda_3(p_{n_G}) \\
\end{bmatrix}.\]
At this time, one can choose $p_j$ as the quadrature points along the boundary $\partial K$. Let $\partial K = e_1 \cup e_2 \cup e_3$. The quadrature points on $e_i$ are denoted by $p_{1,e_i},\cdots, p_{n_g,e_i}$.
In this case, $n_G = 3n_g$, and one can specify the values $\lambda_i(p_j)$ by using the 1-D quadrature points $r_1, r_2, \cdots, r_{n_g}$. The Gaussian quadrature points and weights $r=[r_1,r_2,\cdots,r_{n_g}]$ and $w = [w_1,w_2,\cdots,w_{n_g}]$ on $[0,1]$ are given by \mcode{quadpts1.m}:
\vspace{-0.8cm}
\begin{lstlisting}
[lambda1d,weight1d] = quadpts1(4);  ng = length(weight1d);
[~,id] = sort(lambda1d(:,1));
lambda1d = lambda1d(id,:); weight1d = weight1d(id);
\end{lstlisting}
Here we use \mcode{sort} to guarantee $r_1<r_2<\cdots<r_{n_g}$. Note that
\[\mcode{lambda1d} = \begin{bmatrix} r_1 & r_{n_g} \\
r_2 & r_{n_g-1} \\
\vdots & \vdots \\
r_{n_g} & r_1 \end{bmatrix},\]
and all the row sums are 1, i.e., $r_1 + r_{n_g} = r_2 + r_{n_g-1} = \cdots = 1$. By the definition of the barycentric coordinate functions,  the coordinates $(\lambda_1,\lambda_2, \lambda_3 = 1-\lambda_1-\lambda_2)$ on the three sides of $K$ are:
\begin{itemize}
  \item the 1-th side: $(0,r_{n_g}, r_1), (0, r_{n_g-1}, r_2), \cdots, (0,r_1, r_{n_g})$;
  \item the 2nd side: $(r_1,0,r_{n_g}), (r_2,0, r_{n_g-1}),\cdots,(r_{n_g},0, r_1)$;
  \item the 3rd side: $(r_{n_g}, r_1,0), (r_{n_g-1},r_2,0),\cdots,(r_1, r_{n_g}, 0)$.
\end{itemize}
Therefore the $n_G = 3n_g$ points can be given as
\vspace{-0.8cm}
  \begin{lstlisting}
function [lambdaBd,weightBd] = quadptsBd(order)
%%Gauss quadrature points along the boundary of triangles

% quadrature on [0,1]
[lambda1d,weight1d] = quadpts1(order);  ng = length(weight1d);
[~,id] = sort(lambda1d(:,1));
lambda1d = lambda1d(id,:); weight1d = weight1d(id);
% quadrature on each side
lambdae1 = [zeros(ng,1), lambda1d(:,2), lambda1d(:,1)];
lambdae2 = [lambda1d(:,1), zeros(ng,1), lambda1d(:,2)];
lambdae3 = 1 - lambdae1 - lambdae2;
% quadrature along the boundary of each element
lambdaBd = [lambdae1; lambdae2; lambdae3];
weightBd = repmat(weight1d,1,3);
  \end{lstlisting}

With these quadrature points one can construct the evaluations of the basis functions. The details are omitted. Please refer to \mcode{Base2dBd.m} for implementation details. For the $\mathbb{P}_1$ element, the basis functions are $\phi_i = \lambda_i$, and the elementwise interior evaluations of $u_h$ are given as follows
\vspace{-0.8cm}
  \begin{lstlisting}
elemuhM = zeros(NT,3*ng);
for p = 1:3*ng
    % interpolation at the p-th quadrture point
    for i = 1:length(phi)
        base = phi{i};
        elemuhM(:,p) = elemuhM(:,p) + uh(elem2dof(:,i)).*base(:,p);
    end
end
  \end{lstlisting}
For other types of finite elements, one just needs to modify the basis functions \mcode{base}. For evaluations of the derivatives, for example, $\partial_xu_h$, one needs to change \mcode{base} to the partial derivative $\partial_x$ of the basis functions. The elementwise exterior evaluations of $u_h$ can be computed as
\vspace{-0.8cm}
  \begin{lstlisting}
uhI = zeros(2*NE*ng,1);
uhI(elemQuadM) = elemuhM;
elemuhP = uhI(elemQuadP);
  \end{lstlisting}
We remark that the exterior evaluations on the domain boundary are zero since the corresponding indices in \mcode{elemQuadM} are less than or equal to $n_g \cdot \mcode{NE}$.

The above discussions are summarized in a subroutine:
\vspace{-0.8cm}
  \begin{lstlisting}
function [elemuhM,elemuhP] = elem2edgeInterp(wStr,Th,uh,Vh,quadOrder)
  \end{lstlisting}
One can obtain the evaluations of $u_h$, $\partial_x u_h$ and $\partial_y u_h$ by setting \mcode{wStr} as \mcode{.val}, \mcode{.dx} and \mcode{.dy}, respectively. For convenience, we also return the elementwise unit normal vectors on the quadrature points, denoted by \mcode{elemQuadnx} and \mcode{elemQuadny}, where the first one is for $n_x$ and the second one is for $n_y$. We remark that the normal vectors are simply the unit outer normal vector of the triangles.

\subsubsection{The computation of the jump integral}

With the above preparations, we are able to compute the jump integral
\[
h_e\| [\partial _{n_e}u_h] \|_{0,e}^2
 = h_e\int_e ( [\partial _xu_h]n_x + [\partial _yu_h]n_y )^2 \text{d}s
 = h_e \cdot |e|\sum\limits_{i = 1}^{n_g} w_i ( [\partial _xu_h]n_x + [\partial _yu_h]n_y )^2 (p_{i,e}).
\]

\begin{itemize}
  \item The first step is to compute the elementwise interior and exterior evaluations of $\partial _xu_h$ and $\partial _yu_h$:
  \vspace{-0.8cm}
  \begin{lstlisting}
%% elementwise interior and exterior evaluations at quadrature points
[elemuhxM,elemuhxP,elemnx,elemny] = elem2edgeInterp('.dx',Th,uh,Vh,quadOrder);
[elemuhyM,elemuhyP] = elem2edgeInterp('.dy',Th,uh,Vh,quadOrder);
  \end{lstlisting}
  \item The elementwise jumps are
  \vspace{-0.8cm}
  \begin{lstlisting}
elem2Jumpx = elemuhxM - elemuhxP;
elem2Jumpy = elemuhyM - elemuhyP;
  \end{lstlisting}
  \item The jump integral can be computed by looping of triangle sides:
    \vspace{-0.8cm}
  \begin{lstlisting}
elemJump = zeros(NT,1);
[~,weight1d] = quadpts1(quadOrder);
ng = length(weight1d);
for i = 1:3 % loop of triangle sides
    hei = he(elem2edge(:,i));
    id = (1:ng)+(i-1)*ng;
    cei = hei;
    neix = elemnx(:,id);
    neiy = elemny(:,id);
    Jumpnx = elem2Jumpx(:,id).*neix;
    Jumpny = elem2Jumpy(:,id).*neiy;
    Jumpn = (Jumpnx+Jumpny).^2;
    elemJump = elemJump + cei.*hei.*(Jumpn*weight1d(:));
end
  \end{lstlisting}
  \item We finally get the local error indicators
  \vspace{-0.8cm}
  \begin{lstlisting}
%% Local error indicator
eta = (abs(elemRes) + elemJump).^(1/2);
  \end{lstlisting}
Here we add \mcode{abs} since the third-order quadrature rule has negative weights.
\end{itemize}

\subsection{Numerical example}

Consider the example given in Subsect.~\ref{subsect:implementationFEM}. We employ the D\"{o}rfler marking strategy with parameter $\theta = 0.4$ to select the subset of elements for refinement.

\begin{figure}[!htb]
  \centering
  \subfigure[]{\includegraphics[scale=0.35]{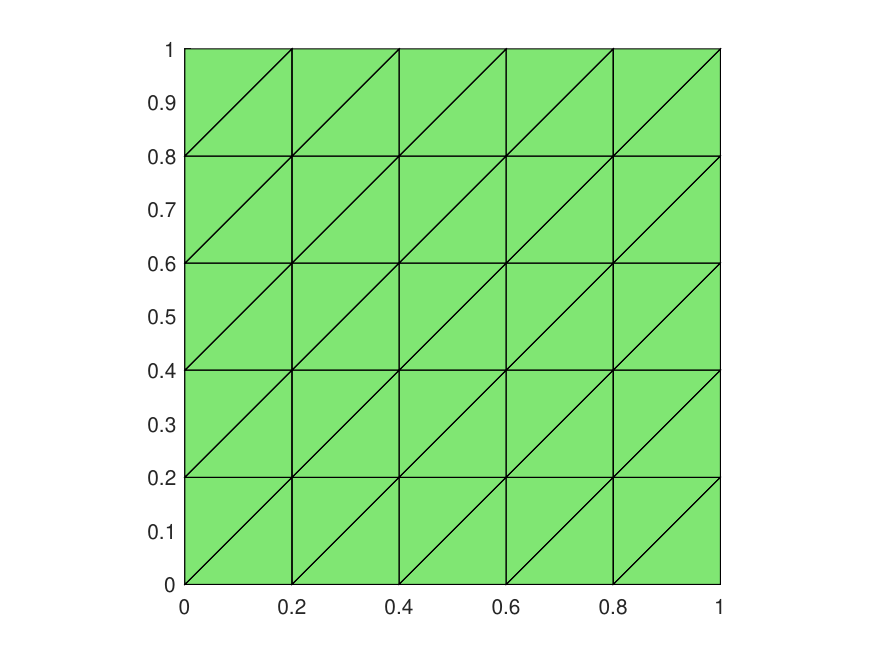}}
  \subfigure[]{\includegraphics[scale=0.35]{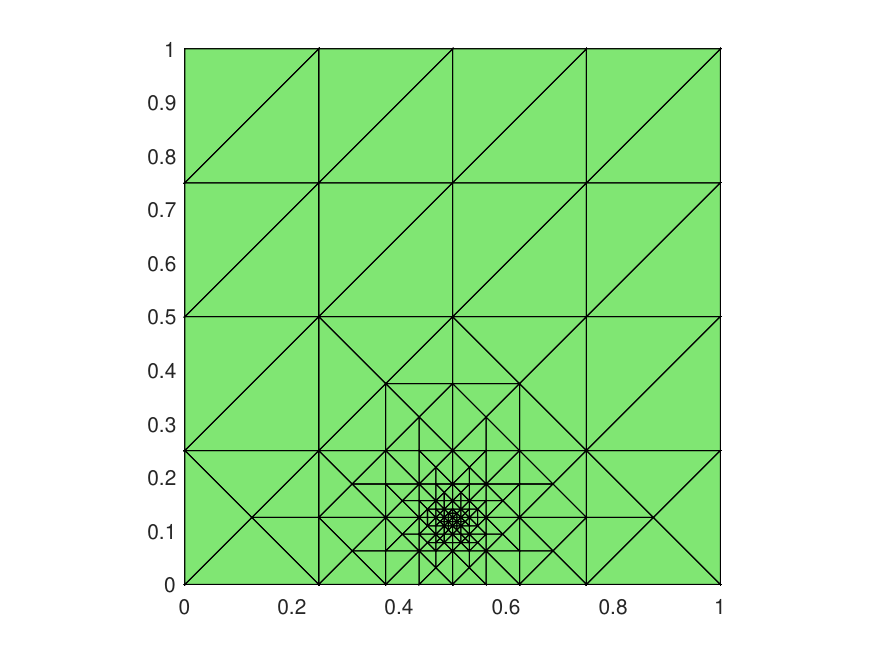}}
  \subfigure[]{\includegraphics[scale=0.35]{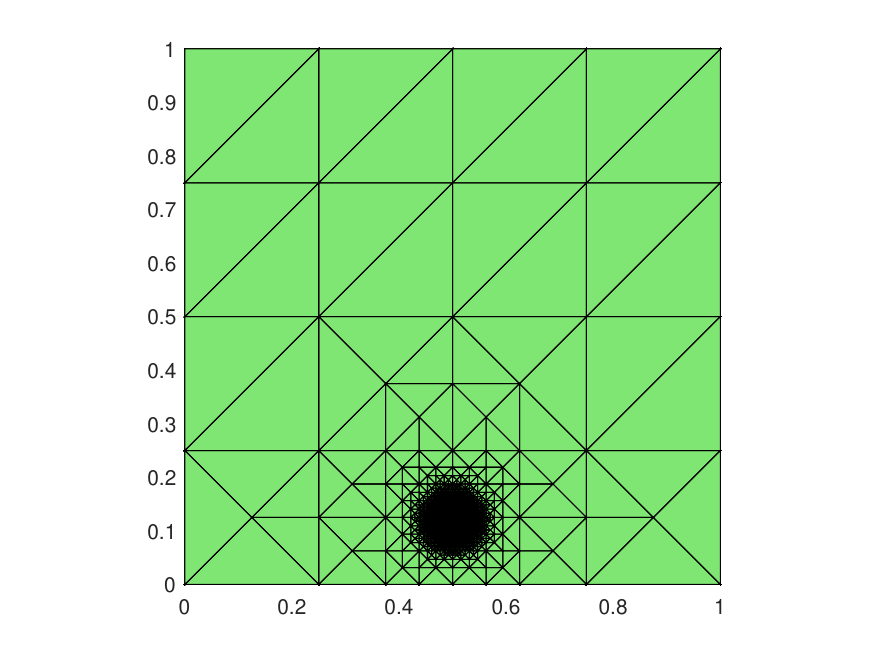}}\\
  \caption{The initial and the final adapted meshes. (a) The initial mesh;
  (b) After 20 refinement steps; (c) After 30 refinement steps}\label{Amesh}
\end{figure}

\begin{figure}[!htb]
  \centering
  \includegraphics[scale=0.45,trim = 40 20 40 20,clip]{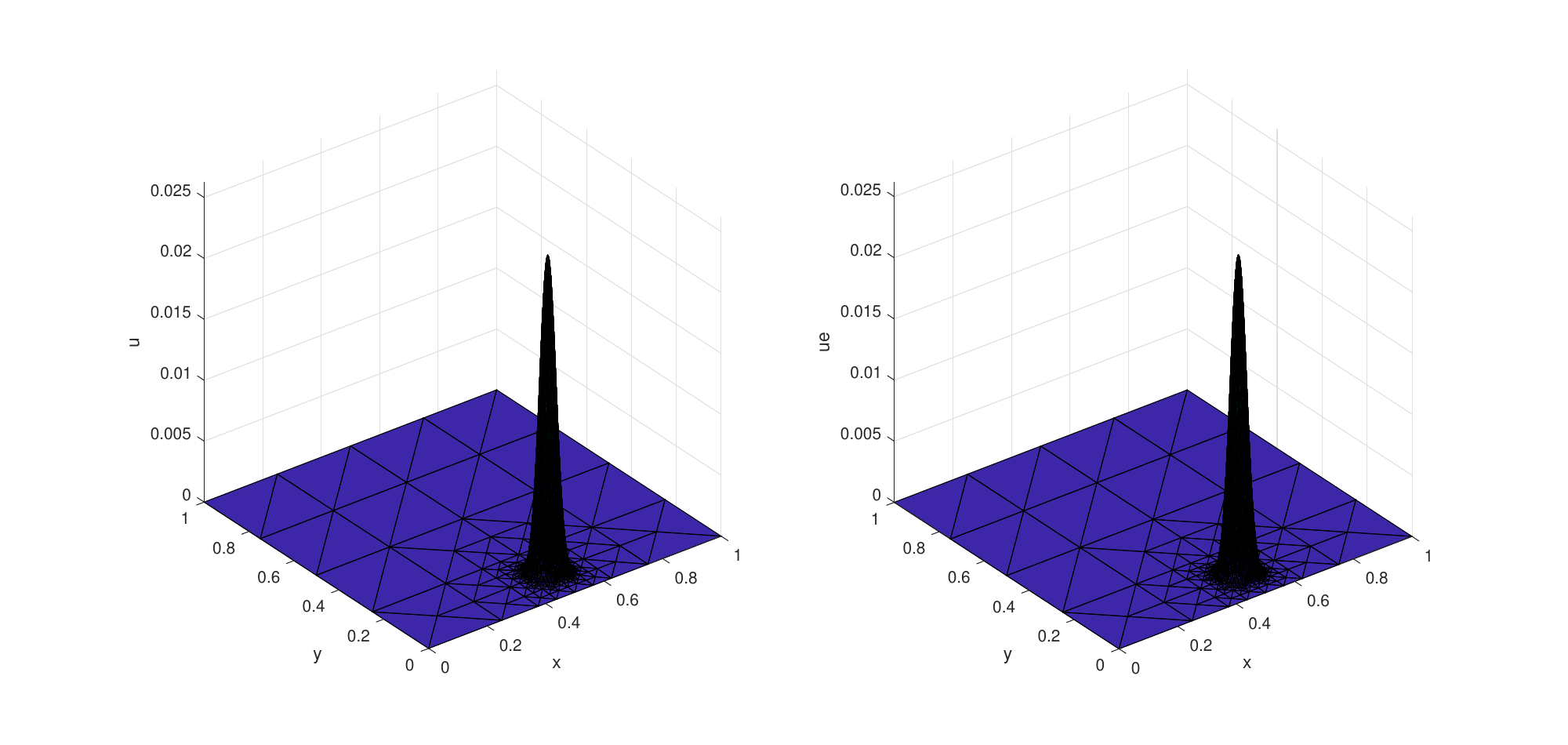}\\
  \caption{The exact and numerical solutions}\label{fig:solutionP1}
\end{figure}

\begin{figure}[!htb]
  \centering
  \subfigure[$\mathbb{P}_1$]{\includegraphics[scale=0.35]{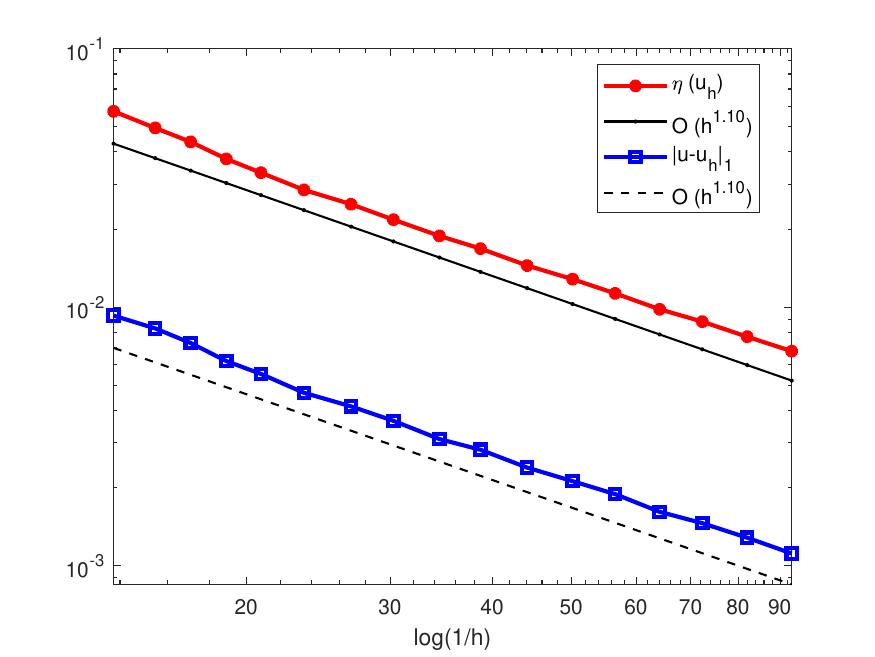}}
  \subfigure[$\mathbb{P}_2$]{\includegraphics[scale=0.35]{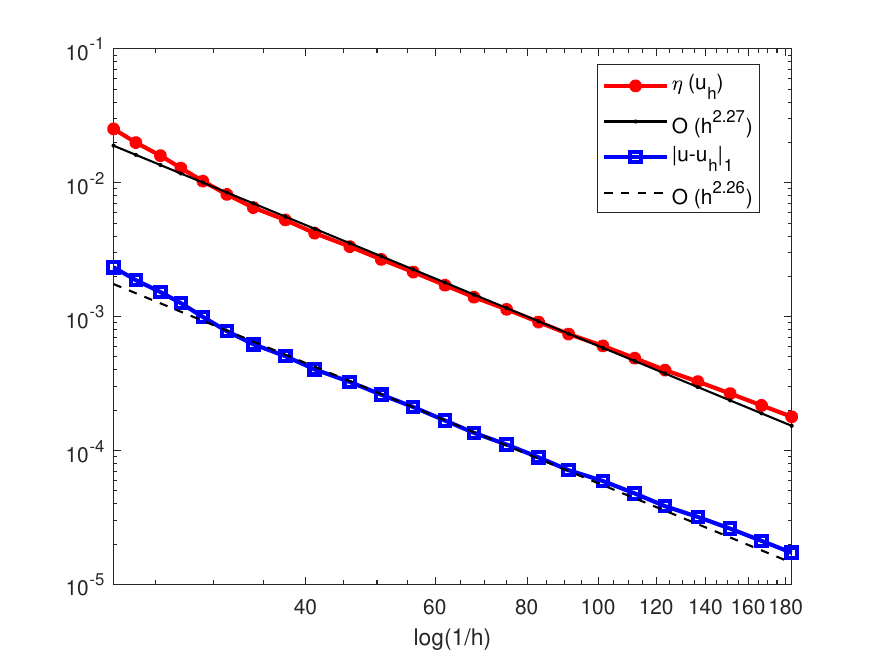}}
  \subfigure[$\mathbb{P}_3$]{\includegraphics[scale=0.35]{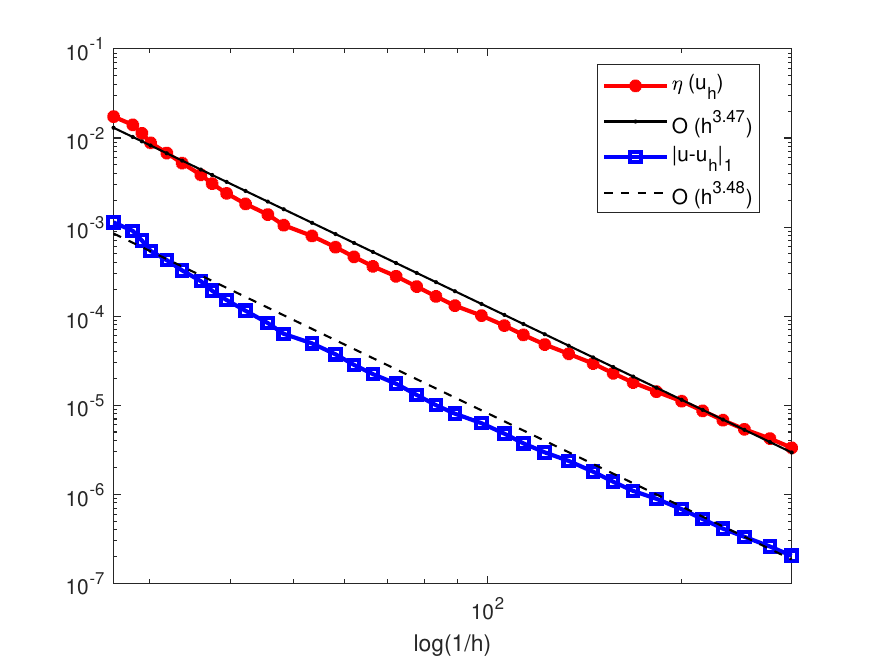}}\\
  \caption{The convergence rates of the error estimators and the errors in $H^1$ norm}\label{AFEMrate}
\end{figure}

The initial mesh and the final adapted meshes after 20 and 30 refinement steps are presented in Fig.~\ref{Amesh}~(a-c),
respectively, where the $\mathbb{P}_1$ element is used. We also plot the adaptive approximation in Fig.~\ref{fig:solutionP1}, which almost coincides with the exact solution. The convergence rates of the error estimators and the errors in $H^1$ norm are shown in Fig.~\ref{AFEMrate}, from which we observe the optimal rates of convergence as predicted by the theory. The full code is available from varFEM package (\url{https://github.com/Terenceyuyue/varFEM}). The subroutine \mcode{PoissonVEM\_indicator.m} is used to compute the local indicator and the test script is \mcode{main\_Poisson\_afem.m}.

\section{$C^0$ interior penalty methods for the biharmonic equation}

\subsection{The introduction to the $C^0$ interior penalty method}

Let $f\in L^2(\Omega)$ and $u\in V:=H_0^2(\Omega)$ be the solution of the biharmonic equation $\Delta^2 u = f$. The variational problem is: Find $u\in V$ such that
\[a(u,v) = \ell(v), \qquad v\in V,\]
where
\[a(u,v) = (\nabla^2 u, \nabla^2 v), \qquad \ell(v) = (f,v). \]
Let $V_h := \mathbb{P}_2(\mathcal{T}_h)\cap H_0^1(\Omega)$ be the quadratic Lagrange finite element space. The symmetric $C^0$ interior penalty method is then to seek $u_h \in V_h$ such that (cf. \cite{Carstensen2022C0IP,Brenner-Sung-2005})
\[a_h(u_h, v_h) = \ell_h(v_h), \qquad v_h \in V_h,\]
where $\ell_h(v_h) = (f, v_h)$ and the bilinear form is
\[a_h(u_h, v_h) := a_{\rm pw}(u_h,v_h) - (\mathcal{J}(u_h,v_h) + \mathcal{J}(v_h,u_h) ) + c_{\rm IP}(u_h,v_h). \]
Here, $a_{\rm pw}$ is the energy scalar product or the interior bilinear form,
\[a_{\rm pw}(u_h,v_h) = \sum\limits_{K\in \mathcal{T}} (\nabla^2 u_h, \nabla^2 v_h)_K.\]
$\mathcal{J}$  is the jump term,
\begin{equation}\label{JumpJ}
\mathcal{J}(u_h,v_h) := \sum\limits_{e\in \mathcal{E}} \int_e\{ \nabla^2 u_h n_e \} \cdot [\nabla v_h] \d s,
\end{equation}
where $\mathcal{E}$ denotes the set of all edges of $\mathcal{T}_h$. The $H^2$ conformity is imposed weakly by the additional penalty term $c_{\rm IP}$,
\[c_{\rm IP}(u_h,v_h) := \sum\limits_{e\in \mathcal{E}}\frac{\sigma_e}{h_e} \int_e [\nabla u_h \cdot n_e] [\nabla v_h \cdot n_e] \d s,\]
where $\sigma_e>0$ is an edge-dependent parameter.

The parameter $\sigma_e$ is used to guaranteed the coercivity or stability of the bilinear form. An automated mesh-dependent selection is proposed in \cite{Carstensen2022C0IP}, given by
\[\sigma_e:= \begin{cases}
\frac{3a h_e^2}{4}\Big( \frac{1}{|K^+|} + \frac{1}{|K^-|} \Big), & \quad e \in \mathcal{E}^0, \\
\frac{3a h_e^2}{|K|} = 2 \cdot \frac{3a h_e^2}{4}\Big( \frac{1}{|K^+|} + \frac{1}{|K^-|} \Big), & \quad e \in \mathcal{E}^\partial,
\end{cases}\]
where we assume $K^+ = K^- = K$ when $e$ is the boundary edge. By Theorem 3.1 there, every choice of $a > 1$ leads to the guaranteed stability.

\begin{remark}
The jump term given in \cite{Brenner-Sung-2005} is
\[\sum\limits_{e\in \mathcal{E}} \int_e \Big\{ \frac{\partial ^2 u_h}{\partial  n_e^2} \Big\}  \Big[\frac{\partial v_h}{\partial n_e} \Big] \d s,\]
which is equivalent to the one in \eqref{JumpJ}. In fact, the $C^0$ continuity leads to $\partial_\tau v_h = 0$, and hence $\nabla v_h  = \partial_n v_h n + \partial_{\tau} v_h \tau  = \partial_n v_h n$, where $n = [n_1,n_2]^T$, $\tau = [t_1,t_2]^T$ and $n \cdot \tau = 0$.
A direct manipulation yields
\begin{align*}
\Big\{ \frac{\partial ^2 u_h}{\partial  n^2} \Big\}  \Big[\frac{\partial v_h}{\partial n} \Big]
& = \{ n^T \nabla^2 u_h n \} [\partial_n v_h] =  ([\partial_n v_h] n)^T \{ \nabla^2 u_h n \} \\
& = [\nabla v_h]^T \{ \nabla^2 u_h n \} = \{ \nabla^2 u_h n \} \cdot [\nabla v_h],
\end{align*}
as required.
\end{remark}

\subsection{The computation of the interior bilinear form}

The shape function space is exactly the $\mathbb{P}_2$-Lagrange finite element space. In contrast to the fully discontinuous Galerkin method, here the degrees of freedom (DoFs) are continuous, which makes the implementation very similar to that for the nonconforming FEMs.

The interior bilinear form
\[a_{\rm pw}(u_h,v_h) = \sum\limits_{K\in \mathcal{T}} (\nabla^2 u_h, \nabla^2 v_h)_K\]
can be easily computed and assembled in varFEM. As for the Poisson equation, we first use \mcode{FeMesh2d.m} to acquire the mesh-related data structures and geometric quantities.
\vspace{-0.8cm}
\begin{lstlisting}
%% Mesh
[node,elem] = squaremesh([0 1 0 1],0.1);  % h1 = h2 = h = 0.1
bdStr = [];
Th = FeMesh2d(node,elem,bdStr);
\end{lstlisting}
The PDE data is defined in \mcode{biharmonicdata\_C0IP.m}. For the current method, the finite element space is specified by \mcode{Vh = 'P2'} as follows.
\vspace{-0.8cm}
\begin{lstlisting}
%% PDE data
pde = biharmonicdata_C0IP();
%% Finite element space
Vh = 'P2'; quadOrder = 5;
\end{lstlisting}
Noting that
\[(\nabla^2 u, \nabla^2 v) = v_{xx} u_{xx} + v_{xy}u_{xy} + v_{yx}u_{yx} + v_{yy}u_{yy},\]
we can compute the stiffness matrix with respect to $a_{\rm pw}$ in the following way:
\vspace{-0.8cm}
\begin{lstlisting}
%% Interior biliner form
Coef  = {1,1,1,1};
Test  = {'v.dxx','v.dxy','v.dyx','v.dyy'};
Trial = {'u.dxx','u.dxy','u.dyx','u.dyy'};
A = assem2d(Th,Coef,Test,Trial,Vh,quadOrder);
\end{lstlisting}

\subsection{The unified implementation of bilinear forms involving jump and average terms}

Different from the jump integral in the a-posteriori error estimator, the jumps and averages of the interior penalty method appear in the bilinear forms. In the sequel we will present a unified design idea for computing and assembling the jump and penalty terms. For clarity, we shall take the jump term
\[\mathcal{J}(u_h,v_h) := \sum\limits_{e\in \mathcal{E}} \int_e\{ \nabla^2 u_h   n_e \} \cdot [\nabla v_h] \d s\]
as an example. It should be pointed out that the above bilinear form is not suitable for elementwise computation because the jump $[\Phi]$ of the global nodal basis function $\Phi$ involves all neighbouring cells. The most natural way is the edge-wise implementation since only the left and right triangles need to be considered.

\subsubsection{The design ideas}

Given an interior edge $e$, the left and right triangles sharing it as an edge will be denoted by $K_1$ and $K_2$, respectively. The union $K_e = K_1 \cup K_2$ is referred to as the macro element of $e$ in what follows. Denote the local basis functions (or their derivatives) on $K_1$ by $\phi_{i_1}, \cdots, \phi_{i_6}$ and on $K_2$ by $\phi_{j_1}, \cdots, \phi_{j_6}$, where the subscripts are the global numbers of the DoFs.
After removing duplicates, the 12 basis functions in fact are associated with only 9 DoFs, which will be ordered from smallest to largest according to the subscripts. The corresponding ``macro basis functions'' are denoted by $\varphi_1, \varphi_2, \cdots, \varphi_9$.

The design idea is based on a simple observation: If the macro element $K_e$ is considered as an element $K$, and the basis function $\varphi$ on $K_e$ is treated as a basis function $\phi$ on $K$, then the jump terms can be computed and assembled in the way as for the interior bilinear form. For this purpose, we can establish the following correspondence:
\begin{align*}
& \mcode{edge2dof} \qquad \longrightarrow \qquad \mcode{elem2dof}, \\
& \mcode{edgeBase} \qquad \longrightarrow \qquad \mcode{Base2d}.
\end{align*}
Specifically, for the jumps and averages on $e$, only the 9 basis functions on the macro elements have contribution to the bilinear forms. For this reason, if $e$ is considered as an element $K$ for the interior bilinear form, and the local basis functions on the macro element are viewed as local basis functions on $K$, then the data structure \mcode{edge2dof} representing the connectivity on the macro element can be established as the data structure \mcode{elem2dof} for the bilinear form of the Poisson equation (see Remark \ref{rem:elem2dof}).
Accordingly, the evaluations of the basis functions can be stored as in \mcode{Base2d.m} (see Remark \ref{rem:Base2d} and the introduction of \mcode{Base2dBd.m} in Subsect.~\ref{subsubsect:MPevaluation}), except that the quadrature points are replaced by the ones on $e$. In addition, for the sake of the vectorized implementation, we assume a virtual exterior triangle for a boundary edge $e$, with all the DoFs indexed by 1 and all the values of the basis functions set to zero, thus not affecting the assembly.

The key point of the program is the construction of \mcode{edgeBase}, which comes down to determining the left and right evaluations of $\varphi_i$ at the quadrature points $z_1, z_2, \cdots, z_{n_g}$.  To do so, let us introduce a matrix, referred to as the {\it macro basis matrix}, given by
\[P = \begin{bmatrix}
\varphi_1(z_1^-)  & \cdots & \varphi_1(z_{n_g}^-)  & \varphi_1(z_1^+)  & \cdots & \varphi_1(z_{n_g}^+) \\
\varphi_2(z_1^-)  & \cdots & \varphi_2(z_{n_g}^-)  & \varphi_2(z_1^+)  & \cdots & \varphi_2(z_{n_g}^+) \\
\vdots            & \vdots & \vdots                & \vdots            & \vdots & \vdots \\
\varphi_9(z_1^-)  & \cdots & \varphi_9(z_{n_g}^-)  & \varphi_9(z_1^+)  & \cdots & \varphi_9(z_{n_g}^+) \\
\end{bmatrix}_{9 \times 2n_g} =: [P^-, P^+],\]
where the first and the second $n_g$ columns store the evaluations on the left and right triangles, respectively. The $\varphi(z^\pm)$ will be obtained from the evaluations of the basis functions $\phi_{i_1}, \cdots, \phi_{i_6}$ and $\phi_{j_1}, \cdots, \phi_{j_6}$. Note that the evaluations of the basis function $\phi$ on the element $K$ are given in counterclockwise order along the boundary $\partial K$.

The idea of constructing the macro basis matrix $P$ is shown in Algorithm \ref{alg:macroP}.
\begin{algorithm}[!htb]
\caption{\textbf{Construction of the macro basis matrix $P$} \label{alg:macroP}}
\begin{enumerate}
  \item Initialization: $P$ is initialized as a zero matrix.
  \item Left and right padding: Let $s_1,\cdots,s_6$ be the local indices of $i_1,\cdots,i_6$ on the macro element. Then
  \[P^-(s_k, :) = [\phi_{i_k}(z_1^-), \cdots, \phi_{i_k}(z_{n_g}^-)], \qquad k = 1,\cdots, 6.  \]
  Similarly, let $t_1,\cdots,t_6$ be the local indices of $j_1,\cdots,j_6$ on the macro element. One has
  \[P^+(t_k, :) = [\phi_{j_k}(z_1^+), \cdots, \phi_{j_k}(z_{n_g}^+)], \qquad k = 1,\cdots, 6.  \]
  \item Reorder: On the left and right triangles, the evaluations of the basis functions are in counterclockwise order. Hence, we need to modify the order of the columns of $P^-$ and $P^+$ according to the orientation of $e$.
\end{enumerate}
\end{algorithm}
Let $\varphi$ be a basis function on the macro element with the local index given by $s$. If $\varphi$ ``belongs to'' the edge $e$, then
\[\varphi(\bb{z}^-) = P^-(s, :), \qquad \varphi(\bb{z}^+) = P^+(r, :), \]
where $\bb{z} = [z_1, z_2, \cdots, z_{n_g}]$. In contrast, if $\varphi$ does not belong to $e$ and suppose it is on the left triangle, then
\[\varphi(\bb{z}^-) = P^-(t, :), \qquad \varphi(\bb{z}^+) = \bb{0}. \]
It is apparent that the initialization and the padding procedure correspond exactly to the above relations.

\subsubsection{The vectorized implementation of the macro basis matrix}

Following the previous idea in Algorithm \ref{alg:macroP}, we are able to construct the macro basis matrix $P$s for any given edge. All these matrices will be collected in the matrix \mcode{edgePhi}, which takes the form of
\[\mcode{edgePhi} = \begin{bmatrix} P_1 \\ P_2 \\ \vdots \\ P_{\text{NE}}\end{bmatrix},
\quad \mbox{$P_k$ is the macro basis matrix of $e_k$}.\]
It is relatively easy to obtain $P_i$ edge by edge, but the large number of loops makes it inefficient. One way of the vectorized implementation is to deal with the $i$-th DoF of all the left or right cells simultaneously. The pseudo-code reads
\vspace{-0.8cm}
\begin{lstlisting}
% loop of basis functions on the left or right elements
for ib = 1:Ndof
    % the current basis functions
    base = phi{ib}; % (NT, 3*ng)
    % indices in edgePhi
    rowL = idMacro(:, ib) + rP; % rP = (0: 9: (NE-1)*9)'
    rowR = idMacro(:, ib+Ndof) + rP;
    % left padding
    edgePhi(rowL, 1:ng) = subMat(base, Lrows, Lcols); % base(k1, (1:ng)+(e1-1)*ng)
    % right padding: k1 ~= k2
    edgePhi(rowR, (1:ng)+ng) = subMat(base, Rrows, Rcols); % base(k2, (1:ng)+(e2-1)*ng)
end
% reorder w.r.t the orientation
edgePhi = edgePhi(ReverseRows, ReverseCols);
\end{lstlisting}
A brief explanation is given as follows:
\begin{itemize}
  \item \mcode{idMacro} determines the local macro indices of the basis functions on the left or right triangles or the local row indices in each $P_k$. Note that $9(k-1)$ should be added if these indices are associated with the $k$-th block of \mcode{edgePhi}, as indicated by the vector \mcode{rP}.
  \item \mcode{e1} records the local side indices of all edges on the left triangle. The matrix \mcode{base} has $3n_g$ columns with the first, the second and the last $n_g$ columns corresponding to the first, the second and the thrid sides of the triangle, where each side records the values at the $n_g$ quadrature points. The column indices corresponding to \mcode{e1} are obviously given by \mcode{(1:ng)+(e1-1)*ng}.
  \item \mcode{base(k1, (1:ng)+(e1-1)*ng)} means to extract the scattered entries in $(i,j)$ locations, where $i$ and $j$ take the values in \mcode{k1} and \mcode{(1:ng)+(e1-1)*ng}, respectively. In Matlab, however, scripts like \mcode{A([1 2], [3 4])} are to extract the submatrix in the intersection of the indexed rows and columns, not the scattered entries. For this purpose, one can convert the $(i,j)$ locations to the linear indices of the implicit column vector. For convenience, we have written a function \mcode{subMat.m} to realize the scattered extraction, where the built-in function \mcode{sub2ind.m} in Matlab is used.
\end{itemize}

With \mcode{edgePhi}, the base matrix given by $\varphi_i$ can be obtained by extracting the $i$-th row of all $P_k$:
\[w_i = \begin{bmatrix} P_1(i,:) \\ P_2(i,:) \\ \vdots \\ P_{\text{NE}}(i,:) \end{bmatrix}, \qquad i = 1,2,\cdots, 9.\]
They are further collected in cell array as
\[\mcode{edgeBaseM} = \{w_1^-, \cdots, w_2^- \}, \qquad \mcode{edgeBaseP} = \{w_1^+, \cdots, w_2^+ \} ,\qquad i = 1,2,\cdots, 9,\]
where
\[w_i^- = w_i(:, ~1:n_g), \qquad w_i^+ = w_i(:, ~(1:n_g)+n_g ).\]
The above discussion will be summarized in the following subroutine
\vspace{-0.8cm}
\begin{lstlisting}
function [edgeBaseM,edgeBaseP,edge2dof] = edgeBase(wStr,Th,Vh,quadOrder)
\end{lstlisting}
To avoid repeated computations, we allow \mcode{wStr} to take more than one string as will be seen in the next subsection.

\subsubsection{The computation and assembly of the jump and penalty terms}

With the previous preparation, stiffness matrices for the jump and penalty terms can be computed and assembled as for the interior bilinear form. Let us take the jump term as an example. A simple calculation shows that
\begin{align*}
& [\nabla v_h] \cdot \{ \nabla^2 u_h   n_e \} \\
& = \begin{bmatrix} [v_{h,x}] \\ [ v_{h,y}] \end{bmatrix} \cdot
\begin{bmatrix} \{ u_{h,xx} n_{e,x} + u_{h,xy} n_{e,y} \}  \\
 \{ u_{h,yx} n_{e,x} + u_{h,yy} n_{e,y} \}
\end{bmatrix}  \\
& = [v_{h,x}] ( \{ u_{h,xx} \}  n_{e,x} + \{ u_{h,xy} \}  n_{e,y} )
+ [v_{h,y}] ( \{ u_{h,yx} \}  n_{e,x} + \{ u_{h,yy} \}  n_{e,y} ) .
\end{align*}

The first step is to obtain the basis matrix. Since the bilinear forms involve the first-order and second-order derivatives, we can compute the left and right evaluations as
\vspace{-0.8cm}
\begin{lstlisting}
%% edgeBase: left and right evaluations of basis functions on macro elements
wStr = {'.dx', '.dy', '.dxx', '.dxy', '.dyy'};
[edgeBaseM,edgeBaseP,edge2dof] = edgeBase(wStr,Th,Vh,quadOrder);
[wxM,wyM,wxxM,wxyM,wyyM] = deal(edgeBaseM{:});
[wxP,wyP,wxxP,wxyP,wyyP] = deal(edgeBaseP{:});
\end{lstlisting}

The sparse assembly index is similar to the one for the interior bilinear form. One just needs to replace \mcode{elem2dof} there by \mcode{edge2dof} (see Remark \ref{rem:elem2dof}).
\vspace{-0.8cm}
\begin{lstlisting}
%% Sparse assembly index
ii = reshape(repmat(edge2dof, Ndofe,1), [], 1);  % Ndofe = 9
jj = repmat(edge2dof(:), Ndofe, 1);
\end{lstlisting}
The \mcode{edge2dof} is of size $\mcode{NE} \times 9$, with each row representing a macro element. For the boundary edges, we have introduced virtual exterior triangles with the DoF indices set to 1. Since the realization is simple, we omit the details. Please refer to the subroutine \mcode{edgeBase.m}.

The stiffness matrices can be computed as follows.
\vspace{-0.8cm}
\begin{lstlisting}
%% Stiffness matrices J and C for jump and penalty terms
KJ = zeros(NE,Ndofe^2);  KC = zeros(NE,Ndofe^2);  % Ndofe = 9
s = 1;
for i = 1:Ndofe
    for j = 1:Ndofe
        % J
        % J: first term
        v1i = wxM{i}-wxP{i};
        u1j = 0.5*Cavg.*(wxxM{j}+wxxP{j}).*nxg + 0.5*Cavg.*(wxyM{j}+wxyP{j}).*nyg;
        % J: second term
        v2i = wyM{i}-wyP{i};
        u2j = 0.5*Cavg.*(wxyM{j}+wxyP{j}).*nxg + 0.5*Cavg.*(wyyM{j}+wyyP{j}).*nyg;
        % J: combined
        KJ(:,s) = he.*sum(ww.*(v1i.*u1j + v2i.*u2j), 2);

        % C
        vi = (wxM{i}-wxP{i}).*nxg + (wyM{i}-wyP{i}).*nyg;
        uj = (wxM{j}-wxP{j}).*nxg + (wyM{j}-wyP{j}).*nyg;
        KC(:,s) = se.*sum(ww.*vi.*uj, 2);

        s = s+1;
    end
end
\end{lstlisting}
Since the exterior values are set as zero and $\{u\} := u$ for the boundary edges, we need to multiply such averages by 2, which is realized by the factor \mcode{Cavg} in the code. For the penalty term, $\sigma_e$ is given by \mcode{se}.

The final stiffness matrix can be assembled by the built-in function \mcode{sparse.m} as
\vspace{-0.8cm}
\begin{lstlisting}
%% Assemble the jump and penalty terms
J = sparse(ii,jj,KJ(:),NNdof,NNdof);
C = sparse(ii,jj,KC(:),NNdof,NNdof);
kk = A - J - J' + C;
\end{lstlisting}

The right-hand side and the boundary conditions can be implemented as for the Poisson equation.
\vspace{-0.8cm}
\begin{lstlisting}
%% Assemble the right-hand side
Coef = pde.f;  Test = 'v.val';
ff = assem2d(Th,Coef,Test,[],Vh,quadOrder);

%% Apply Dirichlet boundary conditions
g_D = pde.g_D;
on = 1;
uh = apply2d(on,Th,kk,ff,Vh,g_D);
\end{lstlisting}

\subsection{Numerical experiment}

\begin{figure}[H]
  \centering
  \subfigure[$h=0.1$]{\includegraphics[scale=0.5]{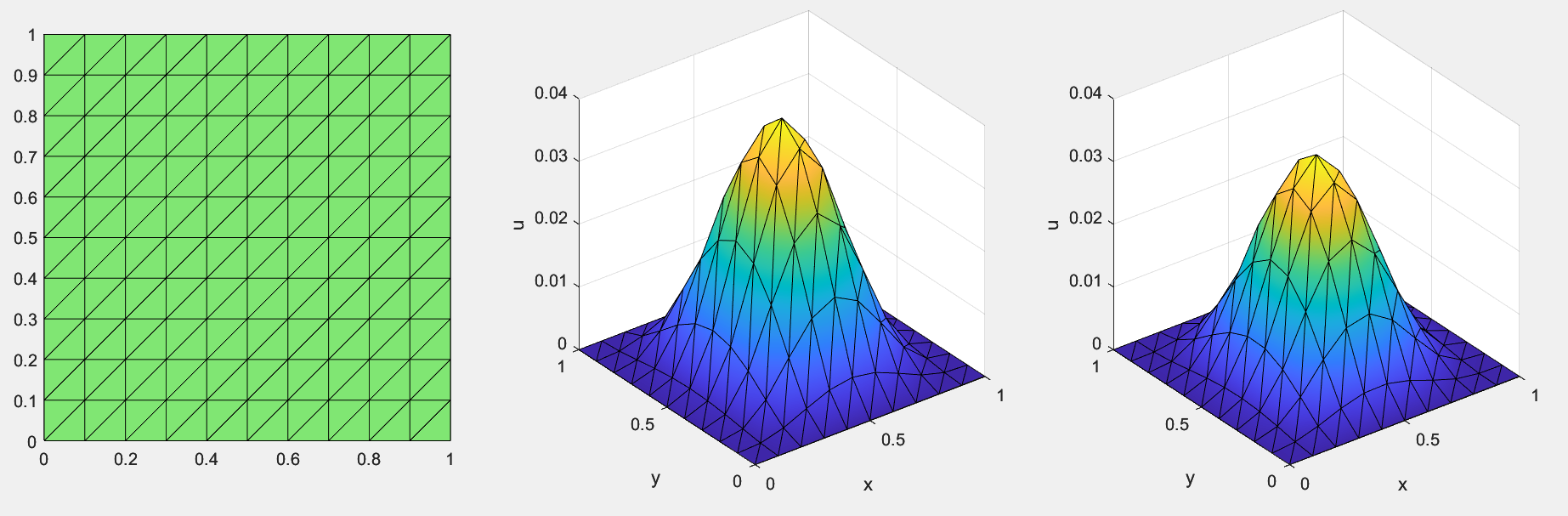}}\\
  \subfigure[$h=0.01$]{\includegraphics[scale=0.5]{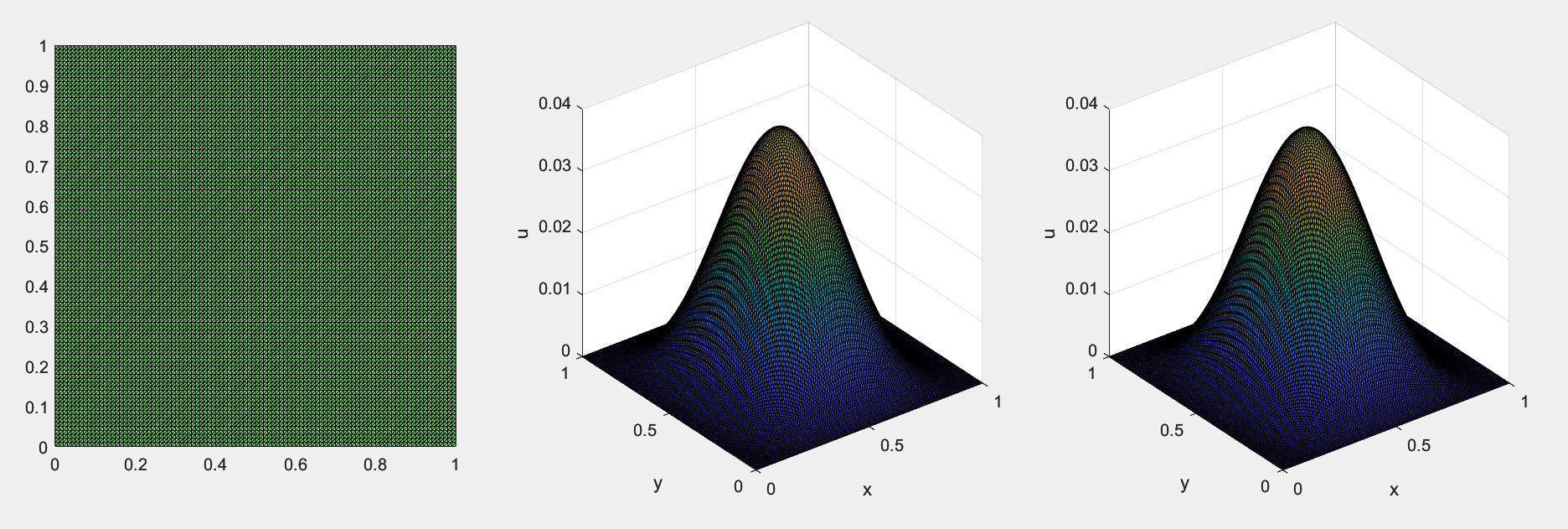}}\\
  \caption{Numerical and exact solutions of the $C^0$ interior penalty method for the biharmonic equation}\label{fig:C0IP}
\end{figure}

Let us consider the biharmonic equation with the exact solution given by $u = 10 x^2 y^2 (1-x)^2 (1-y)^2 \sin(\pi x)$. The numerical solutions for the mesh size $h = 0.1$ and $h = 0.01$ are displayed in Fig.~\ref{fig:C0IP}, which are well matched with the exact solutions and are almost same with the results in \cite{Carstensen2022C0IP}.

We remark that the implementation in \cite{Carstensen2022C0IP} involves some advanced data structures, which makes the code not easy to read and understand. The implementation in this article is more concise and unified for other jumps and averages. Due to the vectorization in constructing the macro basis matrices, our code is also of high efficiency.
For the unit square, when $h=0.01$, the computational times of our implementation and the one in \cite{Carstensen2022C0IP} are about 3.9 seconds and 2.6 seconds, respectively, where the slight overhead is due to the repeated computations in the several modules.

\section{Concluding remarks} \label{sec:conclude}

In this paper, unified implementations of the adaptive FEMs and $C^0$ interior penalty methods were presented. For the adaptive FEMs, the jumps and averages are associated with the known finite element functions, while the $C^0$ interior penalty method for the biharmonic equation is a typical example with jumps and averages in the bilinear forms. The design ideas can be extended to other types of Galerkin methods. The package is accessible on \url{https://github.com/Terenceyuyue/varFEM}.

\end{document}